\documentclass[journal]{IEEEtran}

\usepackage{color}
\usepackage{amsthm}
\usepackage{amsmath}
\usepackage{amssymb}
\usepackage{hyperref}
\usepackage{mathrsfs}
\usepackage{cite}
\usepackage{rgsMacros}
\usepackage{wrapfig}
\allowdisplaybreaks 

\usepackage{hyperref}

\let\labelindent\relax
\usepackage{enumitem}
\usepackage{balance}
\usepackage{enumitem}

\definecolor{olivegreen}{rgb}{0.14,0.29,0}

\newtheorem{exe}{Example}
\newtheorem{prob}{Problem}
\newtheorem{corol}{Corollary}
\newtheorem{ass}{Assumption}
\newtheorem{defin}{Definition}
\newtheorem{cla}{Claim}
\newtheorem{rem}{Remark}
\newtheorem{lem}{Lemma}
\newtheorem{prop}{Proposition}
\newtheorem{thm}{Theorem}
\newtheorem{fct}{Fact}
\newenvironment{lemma}{\begin{lem}}{  \end{lem}}

\newenvironment{corollary}{\begin{corol}}{ \end{corol}}
\newenvironment{example}{\begin{exe}}{ \end{exe}}
\newenvironment{remark}{\begin{rem}}{\hfill \end{rem}}

\newenvironment{assumption}{\begin{ass}}{\hfill  \end{ass}}
\newenvironment{theorem}{\begin{thm}}{ \end{thm}}
\newenvironment{definition}{\begin{defin}}{\hfill 
 \end{defin}}

\usepackage{wrapfig}
\usepackage{chngcntr}
\usepackage{apptools}
\AtAppendix{\counterwithin{lem}{section}}
\AtAppendix{\counterwithin{thm}{section}}
\AtAppendix{\counterwithin{corol}{section}}
\usepackage{mathtools}
\usepackage{comment}
\usepackage{graphicx,xspace,epsfig,syntonly,dsfont}
\usepackage{bbm}
\newcommand\NoIndent[1]{%
  \par\vbox{\parbox[t]{\linewidth}{#1}}%
}

\usepackage[belowskip=-15pt,aboveskip=0pt]{caption}

\usepackage{dsfont}
\newif\ifitsdraft
\def\itsdraft{\global\itsdrafttrue}
\itsdraft  
\begin{document}

\title{\LARGE \bf  Optimal Safety for Constrained Differential Inclusions using Nonsmooth Control Barrier Functions }

\author{
\IEEEauthorblockN{
Masoumeh Ghanbarpour, Axton Isaly,  Ricardo G. Sanfelice, Warren E. Dixon}
\thanks{M. Ghanbarpour and R. G. Sanfelice are with the Department of Electrical and Computer Engineering, University of California, 1156 High Street, Santa Cruz, CA 95064, USA. maghanba@ucsc.edu, ricardo@ucsc.edu.
Research partially supported by the National Science Foundation under Grant no. ECS-1710621, Grant no. CNS-1544396, and Grant no. CNS-2039054, by the Air Force Office of Scientific Research under Grant no. FA9550-19-1-0053, Grant no. FA9550-19-1-0169, and Grant no. FA9550-20-1-0238, and by the Army Research Office under Grant no. W911NF-20-1-0253.

Axton Isaly and Warren E. Dixon are with the Department of Mechanical
and Aerospace Engineering, University of Florida, Gainesville FL 32611-6250, USA. Email: {axtonisaly1013, wdixon}@ufl.edu}
 }

\maketitle
\begin{abstract}
For a broad class of nonlinear systems, we formulate the problem of guaranteeing safety with optimality under constraints.
Specifically, we define controlled safety for differential inclusions with constraints on the states and the inputs. 
Through
the use of nonsmooth analysis tools, we show that a continuous
optimal control law can be selected from a set-valued constraint
capturing the system constraints and conditions guaranteeing
safety using control barrier functions. Our results guarantee
optimality and safety via a continuous state-feedback law designed using nonsmooth control barrier functions.
An example pertaining to obstacle avoidance with a target
illustrates our results and the associated benefits of using nonsmooth control barrier
functions.

\end{abstract}   
\vspace{-5mm}
\section{introduction}  
A powerful approach to guarantee safety for a dynamical system without computing the solutions consists of using barrier functions. Stemming from optimization theory and seminal work by Nagumo in \cite{nagumo1942lage}, a barrier function ensures that, when properly initialized, the solutions of the dynamical system do not reach an unsafe set.
This approach has been exploited for the study of continuous-time, discrete-time, and hybrid systems; see, e.g., \cite{prajna2004safety}; \cite{agrawal2017discrete}; and \cite{kong2013exponential,maghenem2021sufficient}, respectively. 
The extension of the barrier function concept to the case when the system has an input, known as control barrier function (CBF), which is instrumented to synthesize control laws has been pursued in \cite{wieland2007constructive,ames2014control} for continuous-time systems, \cite{agrawal2017discrete} for discrete-time systems, and in \cite{glotfelter2019hybrid}, \cite{isaly2021encouraging} for differential inclusions. To facilitate finding suitable barrier functions needed in control applications, the work in \cite{maghenem2019multiple, glotfelter2017nonsmooth} proposed multiple and nonsmooth barrier functions, \blue{but does not consider optimality}.

Recent developments in combining optimization techniques and safety constraints have led to optimization problems that, when solved numerically, result in a control law that assures both safety and optimality. In \cite{ames2014control}, the authors proposed a
quadratic program to find a minimum norm control law that
ensures safety and stability for nonlinear control affine differential equations. Though powerful, \blue{continuity of} the resulting feedback law is not well characterized. In our recent work, we consider feasibility and continuity of the feedback control law defined by the multiple continuously differentiable barrier functions \cite{isaly2022feasibility}.
Here, we consider differential inclusions with constraints, which are more general than differential equations. 
\blue{Differential inclusions are effective at modeling dynamical systems with uncertainty.}
Nonsmooth barrier functions emerge naturally in many control problems, such as obstacle avoidance. 
The work in \cite{glotfelter2019hybrid} and \cite{glotfelter2017nonsmooth} consider nonsmooth CBFs for differential inclusions, but do not guarantee optimality or continuity of the control law.
Motivated by the need for feedback controllers that assure safety and optimality with good regularity properties, 
\blue{we propose methods to design optimal state-feedback laws using nonsmooth barrier functions that, notably, are continuous.
Specifically,  for constrained differential inclusions, } we propose sufficient conditions for selecting a continuous safe and optimal control law by minimizing a cost function and using a nonsmooth CBF.
More precisely,
we consider differential inclusions with state and input constraints given as
\begin{equation} \label{eq:sys}
   \dot{x} \in F(x,u) \quad (x, u) \in C:= C_x \times C_u 
\end{equation}
where $F:\mathbb{R}^n \times \mathbb{R}^m \rightrightarrows \mathbb{R}^n$, 
 $C_x \subset \mathbb{R}^n$, $ C_u \subset \mathbb{R}^m$.
Let $X_u \subset\mathbb{R}^n$ be unsafe set and $X_o \subset C_x$ indicate the desirable initial set. Using nonsmooth CBFs w.r.t. $(X_o,X_u)$ for constrained differential inclusions as in \eqref{eq:sys}, we provide sufficient conditions for the existence of a continuous safe control law that minimizes a cost function over the set-valued map providing safe inputs for each state. More precisely, the problem we study consists of solving
\begin{equation} \label{eq:opt1}
\begin{aligned}
  \kappa^*(x)=\argmin_{} \quad &  \mathcal{L}(x,u)\\
\textrm{s.t.} \quad  u \in D(x) \\
\end{aligned}
\end{equation}
for each $x \in C_x$ so as to synthesize the optimal state-feedback law $\kappa^*$ ensuring safety, where $\mathcal{L}:\mathbb{R}^n \times \mathbb{R}^m \to \mathbb{R}$ denotes the cost function and the set-valued map $D:\mathbb{R}^n \rightrightarrows \mathbb{R}^m$ indicates feasible safe control inputs \blue{at the current state}.
\ifitsdraft
Our work reveals key properties for the cost function and for the constraints such that the resulting optimal control law is continuous. 
\else 
\fi

\noindent \textbf{Contributions.}
This paper makes the following contributions:
\begin{enumerate}
    \item In Theorem \ref{thm1}, we specify sufficient conditions for the existence of the continuous safe control law using a nonsmooth CBF for \eqref{eq:sys}. 
    \item In Lemma \ref{lm:3} and Lemma \ref{lm:opt_select}, we specify conditions such that the resulting map $\kappa^*$ in \eqref{eq:opt1} is single valued and continuous.
    \item In Theorem \ref{lm:safeOptCont}, we formulate sufficient conditions to obtain a continuous and safe state-feedback law that minimizes the given cost function and meets the safety constraints.
\end{enumerate}
\ifitsdraft 
The remainder of the paper is organized as follows. Preliminary notions are in Section \ref{sec:prelim}. Solutions to differential inclusions are recalled in Section \ref{sec:diff_Inc}. The definitions for safety and CBF are given in Section \ref{sec:safe_N}. Sufficient conditions for
safety using CBF are formulated in Section \ref{sec:suff_con}. In the sections \ref{sec:cont_optimal} and \ref{sec:cont_opt_safe}, we present sufficient conditions for continuous optimal and continuous optimal safe control law, respectively.
\else
\fi
\textbf{Notation.} 
For $x$, $y \in \mathbb{R}^n$, $x^{\top}$ denotes the transpose of $x$, $|x|$ the Euclidean norm of $x$, and $\langle x, y \rangle$ denotes the inner product between $x$ and $y$; namely, $\langle x, y \rangle:= x^\top y$. For a set $K \subset \mathbb{R}^n$, we use $\operatorname{int}(K)$ to denote its interior, $\partial K$ its boundary, $\mbox{cl}(K)$ its closure, $U(K)$ to denote an open neighborhood around $K$, and $\operatorname{co}(K)$ to indicate the convex hull of the set $K$. For $O \subset \mathbb{R}^n$, $K \backslash O$ denotes the subset of elements of $K$ that are not in $O$. For a function $f : \mathbb{R}^n \rightarrow \mathbb{R}^m$, $\dom f$ denotes the domain of definition of $f$, $\operatorname{Graph}(f)$ indicates graph of $f$, and if $f$ is continuously  differentiable, $\nabla f$ denotes the gradient of $f$. If $f$ is locally Lipschitz, $\partial_C f$ denotes Clarke generalized gradient of $f$.
Let $\mathbb{B}$ denote the closed unit ball in $\mathbb{R}^n$ centered at the origin. By $F : \mathbb{R}^m \rightrightarrows \mathbb{R}^n $, we denote a set-valued map associating each element 
$x \in \mathbb{R}^m$ into a subset $F(x) \subset \mathbb{R}^n$. 
\vspace{-3 mm}
\section{Preliminaries}  \label{sec:prelim}
\subsection{Basic Definitions} 
\ifitsdraft
The semicontinuity definitions are from \cite{rockafellar2009variational,goebel2012hybrid,aubin2009set,freeman2008robust}.
\begin{definition} 
(Semicontinuous set-valued maps) \label{def:semi}
Consider a set-valued map $F: C \rightrightarrows \mathbb{R}^n$, where $C \subset \mathbb{R}^m$.
\begin{itemize}
\item The map $F$ is said to be \emph{outer semicontinuous} at $x \in C$, if for every sequence $\left\{x_i\right\}^{\infty}_{i=0} \subset C$ and for every sequence  
$\left\{ y_i \right\}^{\infty}_{i=0} \subset \mathbb{R}^n$ with 
$\lim_{i \rightarrow \infty} x_i = x$, $\lim_{i \rightarrow \infty} y_i = y \in \mathbb{R}^n$, and $y_i \in F(x_i)$ for all $i \in \mathbb{N}$, we have $y \in F(x)$. 
\item The map $F$ is said to be \emph{lower semicontinuous} (or equivalently, \emph{inner semicontinuous}) at $x \in C$ when for every open set $V \subset \text{Im}(F)$ such that $F(x) \cap V \not = \emptyset$, there exists $U(x)$ such that for each $z \in U(x)$, $F(z) \cap V \not = \emptyset$.
\item The map $F$ is said to be \emph{upper semicontinuous} at $x \in C$ if, for each 
$\varepsilon > 0$, there exists $U(x)$ such that for each $y \in U(x) \cap C$, $F(y) \subset F(x) + \varepsilon \mathbb{B}$.
\item The map $F$ is said to be \emph{continuous} at $x \in C$ if it is both outer and lower semicontinuous at $x$.
\end{itemize}

Furthermore, the map $F$ is said to be upper, lower, outer semicontinuous, or continuous, if, respectively, it is upper, lower, outer semicontinuous, or continuous for all $x \in C$.
\end{definition}
\begin{remark}
Based on \cite[Lemma $5.15$]{goebel2012hybrid}, every upper semicontinuous set-valued map with closed values is outer semicontinuous.
Conversely, every outer semicontinuous and locally bounded set-valued map is upper semicontinuous with compact images.
\end{remark}
\begin{definition}(Semicontinuous single-valued maps) \label{def:semi_sing}
Consider a scalar function $f: C \rightarrow \mathbb{R}$, where $C \subset \mathbb{R}^m$.
\begin{itemize}
\item The scalar function $f$ is said to be \emph{lower semicontinuous } at $x \in C$ if, for every sequence $\left\{ x_i \right\}_{i=0}^{\infty} \subset C$ such that $\lim_{i \rightarrow \infty} x_i = x$, we have $\liminf_{i \rightarrow \infty} f(x_i) \geq f(x)$. 
\item The scalar function $f$ is said to be \emph{upper semicontinuous } at $x \in C$, if for every sequence $\left\{ x_i \right\}_{i=0}^{\infty} \subset C$ such that $\lim_{i \rightarrow \infty} x_i = x$, we have $\limsup_{i \rightarrow \infty} f(x_i) \leq f(x)$.  
\item The scalar function $f$ is said to be \emph{continuous} at $x \in C$ if it is both upper and lower semicontinuous at $x$. \end{itemize}
Furthermore, $f$ is said to be upper, lower semicontinuous, or continuous, if respectively, it is upper, lower semicontinuous, or continuous for all $x \in C$.
\end{definition}
\begin{definition} (Clarke generalized gradient \cite[Theorem $2.5.1$]{clarke1990optimization}) \label{def:grad:cl}
Let $B : \mathbb{R}^n \rightarrow \mathbb{R}$ be a locally Lipschitz function. Let $\Omega \subset \mathbb{R}^n$ be any subset of zero measure, and let $\Omega_B \subset \mathbb{R}^n$ be the set of points at which $B$ fails to be differentiable. Then, the \emph{Clarke generalized gradient} of $B$ at $x$ is given by
\begin{align*}
\partial_C B(x) := \co \left\{ \lim_{i \rightarrow \infty} \nabla B(x_i) : x_i \rightarrow x,~x_i \notin \Omega_B,~x_i \notin \Omega \right\}.
\end{align*}
\end{definition}
\else
\fi

\begin{definition} \label{def:levelbounded}
(Sublevel Bounded \cite[Definition $1.8$]{rockafellar2009variational}) A function $f:\mathbb{R}^n \to \mathbb{R}\cup \pm \infty$ is \emph{sublevel bounded} if for every finite $\alpha \in \mathbb{R}$ the set $\{x\in \mathbb{R}^n\,:\, f(x) \leq \alpha\}$ is bounded.
\end{definition}

\begin{definition} \label{def:levelcoercieve}
(Level Coercive \cite[Definition $3.25$]{rockafellar2009variational}) A function $f:\mathbb{R}^n \to \mathbb{R}\cup \pm \infty$ is \emph{level coercive} if it is bounded below on bounded sets and satisfies $\lim_{|x|\to \infty} \inf \frac{f(x)}{|x|}>0$.
\end{definition}
\ifitsdraft
\begin{definition} (Convex set-valued map \cite[Section $3.4$]{freeman2008robust}) \label{def:convex_setvalued}
The set-valued map $F:X\rightrightarrows Z$, where $X$ is convex, is \emph{convex} if for each $\theta \in [0,1]$ and each $x, y \in X$ we have
$F(\theta x + (1-\theta)y) \subset \theta F(x) + (1-\theta)F(y)$
\end{definition}
\else
\fi
\vspace{-5mm}
\subsection{Differential Inclusions} \label{sec:diff_Inc}
Consider differential inclusions with state constraints
\begin{equation} \label{eq:sys_auto}
   \Sigma: \quad \dot{x} \in \widetilde{F} (x) \quad x \in C_x
\end{equation}
where $\widetilde{F} :\mathbb{R}^n \rightrightarrows \mathbb{R}^n$ and $C_x \subset \mathbb{R}^n$.
Next, solutions to the constrained differential inclusion $\Sigma$ in \eqref{eq:sys_auto} are defined. 
\begin{definition} (Concept of Solution to $\Sigma$) \label{def:Sol}
The function $x : \dom x \rightarrow \mathbb{R}^n$, where $\dom x \subset [0,\infty)$, is a solution to $\Sigma$ in \eqref{eq:sys_auto} if 
\ifitsdraft
\begin{enumerate}[label=(\roman*), series = tobecont, itemjoin = \quad]
\item $x(0) \in \cl(C_x)$ 
\item $t \mapsto x(t)$ is locally absolutely continuous,
    \item $x(t) \in C_x$ for all $t \in \operatorname{int}(\dom x)$,
\item $ \dot{x}(t) \in \widetilde{F} (x(t))$ for almost all $t \in \dom x$.
\end{enumerate}
\else
(i) $x(0) \in \cl(C_x)$, (ii) $t \mapsto x(t)$ is locally absolutely continuous, (iii)
 $x(t) \in C_x$ for all $t \in \operatorname{int}(\dom x)$, (iv) $ \dot{x}(t) \in \widetilde{F} (x(t))$ for almost all $t \in \dom x$.
\fi
\end{definition}

A solution $t \mapsto x(t)$ to $\Sigma$ in \eqref{eq:sys_auto} is said to be \emph{complete} if $\dom x$ is unbounded. Furthermore, it is said to be \emph{maximal} if a solution $y$ to $\Sigma$ does not exist such that $x(t) = y(t)$ for all $t \in \dom x$ with $\dom x$ a proper subset of $\dom y$.

Consider the constrained differential inclusion in \eqref{eq:sys}.
Let $\mathcal{U}:C_x \rightrightarrows C_u$ be the set-valued map that provides admissible values in $C_u$ to each $x \in C_x$, namely, $\mathcal{U}(x)$ indicates the feasible control inputs associated with $x$. 
Let $\Sigma_{\mathcal{U}}$ indicate the differential inclusion in \eqref{eq:sys} explicitly constraining $u$ to take values from the set-valued feedback map $\mathcal{U}:C_x \rightrightarrows C_u$, 
\begin{equation} \label{eq:sys_u}
   \Sigma_{\mathcal{U}}: \quad \dot{x} \in F(x,u) \quad u \in \mathcal{U}(x), \,\, x \in C_x.
\end{equation}
\ifitsdraft
Solutions to $\Sigma_{\mathcal{U}}$ are defined as follows.
\begin{definition} (Concept of Solution to $\Sigma_{\mathcal{U}}$) \label{def:Sol_u}
Given a pair of functions $x : \dom x \rightarrow \mathbb{R}^n$ and $u : \dom u \rightarrow \mathbb{R}^m$, where $\dom (x,u) :=\dom x= \dom u \subset [0,\infty)$, $(x,u)$ is a solution to $\Sigma_{\mathcal{U}}$ in \eqref{eq:sys_u} if 
\begin{itemize}
    \item $x(0) \in \cl(C_x)$,
    \item $t \mapsto x(t)$ is locally absolutely continuous,
    \item $t \mapsto u(t)$ is Lebesgue measurable and locally essentially bounded,
    \item $(x(t),u(t)) \in C$ for all $t \in \operatorname{int}(\dom x)$, 
    \item $u(t) \in \mathcal{U}(x(t))$ for all $t \in \dom u$, and
    \item $ \dot{x}(t) \in F(x(t),u(t))$ for almost all $t \in \dom x$.
\end{itemize}
\end{definition}
\else
\fi
Given a selection 
$\kappa(x) \in \mathcal{U}(x)$ for each $x \in C_x$, the (closed-loop) differential inclusion is given by
\begin{equation} \label{eq:sys_selec}
   \Sigma_{\kappa}: \quad \dot{x} \in F(x,\kappa(x)) \quad x \in C_x.
\end{equation}
\vspace{-6mm}
\section{Controlled Safety Notions and CBFs} \label{sec:safe_N}

\subsection{Safety Notions}
Given $\Sigma$ in \eqref{eq:sys_auto} and a set $K \subset C_x$ following \cite{chai2020forward}, we employ the following forward pre-invariance and controlled pre-invariance notions.
 
\begin{definition}
(Forward pre-Invariance) 
The set $K$ is said to be \emph{forward pre-invariant} for the differential inclusion $\Sigma$ in \eqref{eq:sys_auto} if for
each $x_o\in K$, each maximal solution $x$ to $\Sigma$ starting from $x_o$ satisfies $x(t) \in K$ for all $t \in \text{dom}\, x$.
\end{definition}

\begin{definition}
(Controlled pre-Invariance)
A set $K$ is \emph{control pre-invariant} for $\Sigma_{\mathcal{U}}$ in \eqref{eq:sys_u} if there exists a selection 
$\kappa(x) \in \mathcal{U}(x)$ for each $ x \in C_x $
such that $K$ is forward pre-invariant for the resulting differential inclusion $\Sigma_{\kappa}$ in \eqref{eq:sys_selec}.
\end{definition}

Suppose that $X_o\cap X_u = \emptyset$. Safety and controlled safety are defined as follows.
\begin{definition}
(Safety)
Given $\Sigma$ in \eqref{eq:sys_auto}, and $X_o \subset C_x$ and $X_u \subset \mathbb{R}^n$ such that $X_o \cap X_u = \emptyset$, $\Sigma$ is \emph{safe} with respect to $(X_o, X_u)$ if for each solution $x$ to $\Sigma$ starting from 
$x_o \in X_o$, we have $x(t) \in \mathbb{R}^n \backslash X_u$ for all $t \in \dom x$.
\end{definition}
\begin{definition}
(Controlled Safety) Given $\Sigma_\mathcal{U}$ in \eqref{eq:sys_u}, and $X_o \subset \text{dom}\, \mathcal{U}$ and $X_u \subset \mathbb{R}^n$ such that $X_o \cap X_u = \emptyset$, $\Sigma_{\mathcal{U}}$ is \emph{controlled safe} with respect to $(X_o, X_u)$, if there exists a selection 
$\kappa(x) \in \mathcal{U}(x)$ for each $ x \in C_x$
such that the resulting differential inclusion $\Sigma_{\kappa}$ in \eqref{eq:sys_selec} is safe with respect to $(X_o, X_u)$.
\end{definition}

\begin{remark} \label{rm:K}
If we can find a forward pre-invariant (controlled pre-invariant) set $K \subset \mathbb{R}^n$ for $\Sigma$ (respectively, $\Sigma_{\mathcal{U}}$), such that $X_o \subset K$ and $K \cap X_u = \emptyset$, then safety (respectively, controlled safety) is verified with respect to $(X_o, X_u)$.
\end{remark}
\ifitsdraft
 \begin{remark}
 The controlled invarince term is mentioned in \cite{Aubin:1991:VT:120830} and \cite{Blanchini:1999:SPS:2235754.2236030}.
 \end{remark}
 \else \fi
\vspace{-7 mm}
\subsection{CBFs}
Following \cite{wieland2007constructive}, we define control barrier candidates and CBFs.
\begin{definition}
(Barrier Function Candidate) A function $B:\mathbb{R}^n \to \mathbb{R}$ is said to be a barrier function candidate with respect to $(X_o, X_u)$ if  
\begin{equation}
\label{eq:candid}
\begin{aligned} 
 B(x)  > 0 \qquad \forall x \in X_u , \qquad 
 B(x)  \leq 0 \qquad \forall x \in  X_o.  
\end{aligned}
\end{equation}
\end{definition}
The zero sublevel set of $B$ is defined as
$K_e:= \{x \in \mathbb{R}^n\,:\, B(x) \leq 0\}.$
Consider $\Sigma_\mathcal{U}$ in \eqref{eq:sys_u}. Suppose $B$ is a barrier function candidate with respect to $(X_o, X_u)$. Let the set $K$ be defined as 
\begin{equation} \label{eq:k}
    K:= K_e \cap C_x.
\end{equation}
If $K$ is controlled pre-invariant for $\Sigma_\mathcal{U}$, then according to Remark \ref{rm:K}, $\Sigma_\mathcal{U}$ is controlled safe. 

\begin{definition} \label{def_cbf}
(Control Barrier Function) A locally Lipschitz barrier function candidate $B$ with respect to $(X_o, X_u)$ is a \emph{CBF } for $\Sigma_{\mathcal{U}}$ in \eqref{eq:sys_u} if there exists a neighborhood of the boundary of the set $K$, $U(\partial K)$, where $K$ is defined in \eqref{eq:k}, such that the following condition holds:
\begin{equation} \label{eq:cbf_llip} 
\begin{aligned}
    \inf_{u \in \mathcal{U}(x)} \sup_{\eta \in F(x,u), \zeta \in \partial_C B(x)} &\langle \zeta, \eta\rangle \leq 0
    \\
    & \forall x \in \left(U(\partial K) \backslash \operatorname{int}(K) \right) \cap C_x.
\end{aligned}
\end{equation} 
\end{definition}

\ifitsdraft
\begin{remark} \label{re:cbf_c1}
If the barrier function candidate is \emph{continuously differentiable}, condition \eqref{eq:cbf_llip} is written as
\begin{equation}  \label{eq:cbf_c1}
\begin{aligned}
    \inf_{u \in \mathcal{U}(x)} \sup_{\eta \in F(x,u)} \langle  \nabla B(x), \eta\rangle & \leq 0 
    \\
    & \forall x \in \left(U(\partial K) \backslash K \right) \cap C_x.
    \end{aligned}
\end{equation}
\end{remark}
\else
\fi
\normalsize
Let the function $B$ be a CBF for $\Sigma_{\mathcal{U}}$ in \eqref{eq:sys_u} with respect to $(X_o, X_u)$, 
defining the set $K$ in \eqref{eq:k}. Suppose $B$ is locally Lipschitz. We define the function $g:C_x \times C_u \to \mathbb{R}$ as
\begin{equation} \label{eq:g_ll}
g(x,u) := \sup \{ \langle  \zeta, \eta\rangle \,:\, \zeta \in \partial_C B(x), \, \eta \in F(x,u) \}.
\end{equation}
For some continuous function $\gamma:C_x \to \mathbb{R}$, we define the set-valued map $D_{\gamma}:C_x \rightrightarrows C_u$ as 
\begin{equation} \label{eq:dgamma}
D_{\gamma}(x) := \{u \in \mathcal{U}(x)\,:\, g(x,u) +\gamma(x) < 0 \}.
\end{equation}
The set-valued map $D$ for each $x$ gives the set of all feasible control inputs, that makes the function $g+\gamma$ negative; therefore, it provides the possible control inputs for which the CBF decreases along the solutions.

\begin{assumption} \label{ass:C}
The set $C=C_x \times C_u$ is closed.
\end{assumption}

\begin{assumption} \label{ass:F} 
The map $F: C \rightrightarrows \mathbb{R}^n $ is upper semicontinuous and $F(x,u)$ is nonempty, compact, and convex for all $(x,u) \in C$.
\end{assumption}
Assumptions \ref{ass:C} and \ref{ass:F} are known as tight requirements in the literature for the existence of solutions and the structural properties for the set of solutions of differential inclusion; see \cite{aubin2012differential, Aubin:1991:VT:120830, clarke2008nonsmooth}. For Assumption \ref{ass:f_convex_u}, \ifitsdraft
(see Definition \ref{def:convex_setvalued}).
\else
see \cite[Section $3.4$]{freeman2008robust}.
\fi 
\begin{assumption} \label{ass:U}
The feedback map $\mathcal{U}:C_x \rightrightarrows C_u$ is lower semicontinuous with nonempty, closed, and convex values. 
\end{assumption}
\begin{assumption} \label{ass:f_convex_u}
The map $F$ is convex in $u$.
\end{assumption}

Based on Michael's Theorem \cite[Theorem $2.18$]{freeman2008robust}, Assumption \ref{ass:U} guarantees the essential requirement to find a continuous selection from $\mathcal{U}$.
In the following lemma, we indicate the regularity of the maps $g$ and $D_\gamma$ under the said assumptions.

\begin{lemma} \label{lm:dgamma2}
Consider $\Sigma_{\mathcal{U}}$ in \eqref{eq:sys_u} such that Assumptions \ref{ass:C}-\ref{ass:U} hold.
Let $B:\mathbb{R}^n \to \mathbb{R}$  be locally Lipschitz and $\gamma:C_x \to \mathbb{R}$ be a continuous function. Suppose the set-valued map $D_{\gamma}:C_x \rightrightarrows C_u$ is defined in \eqref{eq:dgamma} and $g:C_x \times C_u \to \mathbb{R}$ is defined in \eqref{eq:g_ll}. Then, the following hold:
\begin{enumerate} [label={\ref{lm:dgamma2}.\alph*},leftmargin=*]
\item \label{item:g_usc} The function $g$ is upper semicontinuous,
\item \label{item:d_lsc}The map $D_{\gamma}$ is lower semicontinuous,
\item \label{item:g_u_lsc} When Assumption \ref{ass:f_convex_u} holds, for each $x \in C_x$ the function $u \mapsto g(x,u)$ is convex and lower semicontinuous.
\end{enumerate}
\end{lemma}

\begin{proof}
 To prove \ref{item:g_usc}, using Assumption \ref{ass:F}, we obtain that $F$ has nonempty compact values, and using \cite[Proposition $2.6.2$]{clarke1990optimization}, we conclude that $\partial_C B$ is nonempty and has convex compact values. Therefore, $g$ is well-defined.
Since $B$ is locally Lipschitz, from \cite[Lemma $4.6$]{sanfelice2007invariance} we conclude that $g$ is upper semicontinuous. 

To prove \ref{item:d_lsc},
from \ref{item:g_usc}, we have that $g$ is upper semicontinuous and since $\gamma$ is continuous, we conclude that $g+\gamma$ is upper semicontinuous. 
Since, by Assumption \ref{ass:U}, $\mathcal{U}$ is a lower semicontinuous, \cite[Corollary 2.13]{freeman2008robust} implies that $D_{\gamma}$ is lower semicontinuous.

To prove \ref{item:g_u_lsc}, first we show that $g$ is convex in $u$.
Using Assumption \ref{ass:f_convex_u}, for each $\theta \in [0,1]$, $x \in C_x$, and $u_1, u_2 \in C_u$, we have
\ifitsdraft
\begin{equation}
\begin{aligned}
  g&(  x , \theta u_1 +  (1-\theta) u_2) \\
  &= 
  \sup_{\zeta \in \partial_C B(x), \, \eta \in F(x,\theta u_1 + (1-\theta) u_2)} \langle  \zeta, \eta\rangle 
   \\
  &\leq \sup_{\zeta \in \partial_C B(x), \,\eta \in \theta F(x, u_1 ) + (1-\theta)F(x,u_2)}  \langle \zeta, \eta\rangle 
  \\
  & = \sup_{\zeta \in \partial_C B(x), \,\eta_1 \in  F(x, u_1 ), \eta_2 \in F(x,u_2)}  \langle \zeta, \theta\eta_1 + (1-\theta) \eta_2\rangle 
  \\
  & \leq  \theta \sup_{\zeta \in \partial_C B(x), \,\eta_1 \in  F(x, u_1 )}  \langle \zeta, \eta_1 \rangle \\
  & \quad + (1-\theta)
  \sup_{ \zeta \in \partial_C B(x), \,\eta_2 \in F(x,u_2)} \langle \zeta , \eta_2 \rangle \\
  & = \theta g(x,u_1) + (1-\theta) g(x,u_2)
\end{aligned}
\end{equation}
\else
\small{\begin{equation*}
\begin{aligned}
  g(  x , \theta u_1 +  (1-\theta) u_2) 
   &\leq  \theta \sup_{\zeta \in \partial_C B(x), \,\eta_1 \in  F(x, u_1 )}  \langle \zeta, \eta_1 \rangle \\
  & + (1-\theta)
  \sup_{ \zeta \in \partial_C B(x), \,\eta_2 \in F(x,u_2)} \langle \zeta , \eta_2 \rangle \\
  &  \leq \theta g(x,u_1) + (1-\theta) g(x,u_2)
\end{aligned}
\end{equation*}}
\normalsize
\fi
\normalsize
To prove that $g$ is lower semicontinuous, let $h:\mathbb{R}^n \times \mathbb{R}^m \times \mathbb{R}^n \to \mathbb{R} $ be defined as $h(x,u,\zeta):= \sup_{ \eta \in F(x,u) } \langle \eta, \zeta \rangle$.
For fixed $\zeta$, the map $\eta \mapsto \langle \eta, \zeta \rangle$ is continuous and convex.
Furthermore, since $F$ has bounded values and $\partial_C B$ is bounded, we conclude that $h$ is bounded. Then, \cite[Theorem $9.4$]{rockafellar1970convex}
 implies that $u \mapsto h(x,u,\zeta)$ is lower semicontinuous. 
 \ifitsdraft
 Therefore, since $g$ is equal to
\begin{align}
    g(x,u) = \sup_{ \zeta \in \partial_C B(x)} h(x,u,\zeta)
\end{align}
\else
Therefore, since
$g(x,u) = \sup_{ \zeta \in \partial_C B(x)} h(x,u,\zeta)$,
\fi
$u \mapsto g(x,u)$ is lower semicontinuous for each $x \in C_x$.
\end{proof}

\begin{remark} \label{rm:ncon_D}
In general, when $B$ is not continuously differentiable, $D_\gamma$ is not necessarily outer semicontinuous, therefore, $D_\gamma$ is not continuous.
Since the Clarke generalized gradient of the locally Lipschitz function $B$, $\partial_C B$, is upper semicontinuous, with correct regularity of $F$, the function $g$ can only be upper semicontinuous.
In general, $D_\gamma$ cannot be continuous without continuity of $g$.
\end{remark}
\vspace{-5 mm}
\section{Sufficient Conditions for Safety} \label{sec:suff_con}
\ifitsdraft
Under \eqref{eq:cbf_llip}, using the CBF in Definition \ref{def_cbf}, we define a set-valued map that indicates feasible and safe feedback control inputs for each $x$. Then, we provide sufficient conditions to guarantee the existence of a continuous selection from this set-valued map to ensure safety for the system.
\else
\fi

\begin{theorem} \label{thm1}
Consider $\Sigma_{\mathcal{U}}$ in \eqref{eq:sys_u} such that Assumptions \ref{ass:C}-\ref{ass:f_convex_u} hold.
Let $B : \mathbb{R}^n \rightarrow \mathbb{R}$ be a locally Lipschitz CBF with respect to $(X_o, X_u) \subset \mathbb{R}^n \times \mathbb{R}^n$ defining the set $K$ in \eqref{eq:k}.
Suppose there exists a neighborhood of $\partial K$, denoted by $U(\partial K) $, such that the set-valued map $D_\gamma:C_x \rightrightarrows C_u$  defined in \eqref{eq:dgamma} for $\gamma$ identically zero is nonempty on $\big( U(\partial K) \backslash \text{int(K)}\big) \cap C_x$.
 Then, there exists a continuous control law $\kappa:C_x \to C_u$ that makes $\Sigma_{\mathcal{U}}$ controlled safe with respect to $(X_o,X_u)$.
\end{theorem}
\begin{proof}
Let $g:C_x \times C_u \to \mathbb{R}$ be defined in \eqref{eq:g_ll}, and from \eqref{eq:dgamma}, let the map $D_0: C_x \rightrightarrows C_u$ for $x \mapsto \gamma(x)=0$ be defined as
\begin{equation} \label{eq:D-zero}
   D_{0}(x) := \{u \in \mathcal{U}(x)\,:\, g(x,u) < 0 \}. 
\end{equation}
Since $B$ is a CBF for $\Sigma_{\mathcal{U}}$, suppose $U_1(\partial K)$ is a neighborhood such that \eqref{eq:cbf_llip} holds. Let $U_2(\partial K)$ be a neighborhood such that $U_2(\partial K) \subset U_1(\partial K)$ and $U_2(\partial K) \subset U(\partial K)$.
From \ref{item:g_usc} and \ref{item:g_u_lsc} in Lemma \ref{lm:dgamma2}, we have that $g$ is continuous in $u$;  therefore, for each $x$, its sublevel sets are closed. 
Then, since also for each $x\in C_x$, $\mathcal{U}(x)$ has closed values, we obtain that $\cl(D_0(x)) = \{u \in \mathcal{U}(x)\,:\, g(x,u) \leq 0 \}=:\bar{D}_0$. 
Let $S = \cl( \big( U_2(\partial K ) \backslash K \big) \cap C_x )$.
To prove that there exists a continuous selection from $ \bar{D}_0$ on $S$, based on \cite[Theorem $2.18$]{freeman2008robust}, we show that $\bar{D}_0$ is lower semicontinuous and nonempty on $S$, with closed convex values.
Since Assumptions \ref{ass:C}-\ref{ass:U} hold, \ref{item:d_lsc} in Lemma \ref{lm:dgamma2} implies that 
$D_{0}$ is lower semicontinuous. Using \cite[Proposition $2.3$]{michael1956continuous}, we conclude that $\bar{D}_0$ is lower semicontinuous.
Furthermore,  $\bar{D}_0$ is nonempty on $S$ and has closed values.
Finally, using Assumption \ref{ass:U}, we obtain that $x \mapsto \mathcal{U}(x)$ has convex values and
using the assumption that $u \mapsto g(x,u)$ is convex, we conclude that $\bar{D}_0$ has convex values. 
Therefore, \cite[Theorem $2.18$]{freeman2008robust}
 implies that there exists a continuous selection $\kappa_1:S \to C_u$ such that
$\kappa_1(x) \in \bar{D}_0(x)$ for each $ x \in S$.
Furthermore, since for each $x\in S$, we have 
$\bar{D}_0(x) \subset \mathcal{U}(x)$,
then $\kappa_1$ is also a selection from $\mathcal{U}$. Then, using Assumption \ref{ass:U} and 
\ifitsdraft
Lemma \ref{lm:michael_extend} in the Appendix, 
\else 
\cite[Proposition $1.4$]{michael1956continuous},
\fi
$\kappa_1$ can be extended continuously to the entire $C_x$.
Let $\kappa:C_x \to C_u$ be the extension of $\kappa_1$.
 Finally, to prove that $\Sigma_{\kappa}$ as defined in \eqref{eq:sys_selec} is safe with respect to $(X_o,X_u)$, from \eqref{eq:cbf_llip} we have 
 \vspace{-2mm}
$$\sup_{\eta \in F(x,\kappa(x)), \zeta \in \partial_C B(x)} \langle \zeta, \eta\rangle \leq 0 \quad \forall x \in U_2(\partial K) \backslash K.$$
Therefore,
$\langle \zeta, \eta\rangle \leq 0$ for each $x \in U_2(\partial K) \backslash K, \, \zeta \in \partial_C B(x)$, and for each $\eta \in F(x,\kappa(x))$.
Hence, using Assumption \ref{ass:F} and the fact that $\kappa$ is continuous, we conclude that $x \mapsto F(x, \kappa(x))$ is upper semicontinuous with nonempty, compact, and convex values.
\ifitsdraft
Using Lemma \ref{lm:appOSC_usc} in the Appendix, 
\else 
Using \cite[Lemma $5.15$]{goebel2012hybrid},
\fi 
we conclude that upper semicontinuous maps with compact images are outer semicontinuous and locally bounded set-valued maps. 
Therefore, 
based on \cite[Theorem $4$]{maghenem2021sufficient}, $K$ is forward pre-invariant for $\Sigma_\kappa$. Then, 
$\Sigma_{\kappa}$ is safe with respect to $(X_o, X_u)$. Therefore, $\Sigma_\mathcal{U}$ is controlled safe with respect to $(X_o, X_u)$. 
\end{proof}
\begin{remark}
For a continuously differentiable function $B$ and single-valued function $f:C_x \times C_u \to \mathbb{R}^n$, the condition
\begin{equation} \label{eq:cond_lit}
    \inf_{u \in \mathcal{U}} \langle \nabla B(x),f(x,u)\rangle \leq -\alpha(B(x))
\end{equation}
has been used in the literature; see, e.g., \cite{glotfelter2019hybrid} and \cite{xu2015robustness}, where $\alpha:\mathbb{R} \to \mathbb{R}$ is an extended class $\mathcal{K}_\infty$ function; namely, $\alpha(0)=0$ and $\alpha$ is strictly increasing.
Note that condition \eqref{eq:cbf_llip}
\ifitsdraft
or \eqref{eq:cbf_c1} 
\else 
\fi
is more general than \eqref{eq:cond_lit}, in the sense that the inequality in \eqref{eq:cbf_llip}
\ifitsdraft
or \eqref{eq:cbf_c1} 
\else 
\fi
does not need to hold on the entire set $C_x$; however, in \eqref{eq:cond_lit}, the safety constraint is imposed globally, though it may get relaxed in the interior of the safe set. 
  Another advantage of \eqref{eq:cbf_llip} is that the barrier function needs only to be locally Lipschitz in a neighborhood of the boundary of $K$.
  \blue{Furthermore, using the presented framework, a combination of multiple intersecting and non-intersecting barriers can be addressed independently.}
\end{remark}
\vspace{-5 mm}
\section{Guaranteeing Continuity of \\the Optimal Solution} \label{sec:cont_optimal}
Given a set-valued map $D:\mathbb{R}^n \rightrightarrows \mathbb{R}^m$ indicating all safe and feasible control actions for each $x$, and the desired cost function $\mathcal{L}$, an optimal control law is given by solving 
\ifitsdraft
the following optimization problem:
\begin{equation} \label{eq:opt1}
\begin{aligned}
  \kappa^*(x)= &\argmin_{}  \quad   \mathcal{L}(x,u)\\
&\quad \textrm{s.t.} \quad  u \in D(x). \\
\end{aligned}
\end{equation}
\else
the optimization problem in \eqref{eq:opt1}.
\fi
In the following, we give two sets of conditions concerning the cost function $\mathcal{L}$ and the set-valued map $D$, such that the optimal control law $\kappa^*$ is continuous.

 Berge's Maximum Theorem \cite[Maximum Theorem]{berge1997topological} provides conditions such that the optimal solution map $\kappa^*$ in \eqref{eq:opt1} is upper semicontinuous and has compact values. In the following lemma, we specify conditions such that the optimal solution to \eqref{eq:opt1} is single valued and continuous.
\begin{lemma} \label{lm:3}
Suppose the function $\mathcal{L}:\mathbb{R}^n \times \mathbb{R}^m \to \mathbb{R}$ is continuous and strictly convex in its second argument.
Let the set-valued map $D:C_x \rightrightarrows C_u$ be continuous and have nonempty and compact convex values.
Then, the function $\kappa^*:C_x \to C_u$, defined in \eqref{eq:opt1}, for each $x \in C_x$ is single valued and continuous.
\end{lemma}
\begin{proof}
Since $\mathcal{L}$ is continuous and the set-valued map $D$ is continuous, and has nonempty compact values, \cite[Maximum Theorem]{berge1997topological} implies that the set-valued map $\kappa^*$ in \eqref{eq:opt1} is nonempty and upper semicontinuous, with compact values. 
Since $\mathcal{L}$ is strictly convex in $u$, \cite[Theorem $2.6$]{rockafellar2009variational} implies that for each $x$, $\kappa^*(x)$ in \eqref{eq:opt1} has at most one value.
Therefore, $\kappa^*$ is nonempty and single valued.
Using 
\ifitsdraft
Lemma \ref{lm:appOSC_usc} in the Appendix, 
\else 
\cite[Lemma $5.15$]{goebel2012hybrid},
\fi
we conclude that upper semicontinuous maps with compact images are outer semicontinuous with locally bounded values. From
\ifitsdraft
Lemma \ref{lm:mul_sing} in the Appendix, 
\else 
\cite[Corollary $5.20$]{rockafellar2009variational},
\fi
we conclude that $\kappa^*$ is continuous.
\end{proof}
\vspace{-2mm}
The following lemma is more general than Lemma \ref{lm:3} in the sense that in Lemma \ref{lm:opt_select}, $\cal L$ can be sublevel bounded in $u$ when $D$ does not have bounded values.
\begin{lemma} \label{lm:opt_select}
Suppose the function $\mathcal{L}:\mathbb{R}^n \times \mathbb{R}^m \to \mathbb{R}$ is proper, lower semicontinuous, convex, and strictly convex in its second argument. 
Let the set-valued map $D:C_x \rightrightarrows C_u$ be continuous, nonempty and have convex values. 
Let the function $\kappa^*:C_x \to C_u$ be defined in \eqref{eq:opt1} for each $x \in C_x$. If one of the following conditions holds,
\begin{enumerate}
    \item $\mathcal{L}$ is sublevel bounded in $u$
    (see Definition \ref{def:levelbounded});
    \item $D$ has bounded values, namely, for each $x \in C_x$, $D(x)$ is bounded.
\end{enumerate}
then $\kappa^*$ is single valued on $C_x$ and continuous on $\operatorname{int}(C_x)$.
\end{lemma}
\begin{proof}
To prove the lemma, we use
\ifitsdraft
Corollary \ref{cor:num1} in the Appendix. 
\else 
\cite[Corollary $7.43$]{rockafellar2009variational}.
\fi
\normalsize
Let the function $\tilde{\delta}:\mathbb{R}^n \times \mathbb{R}^m \to \bar{\mathbb{R}}$ be defined as
$$\tilde{\delta}(x,u):= \delta_{D(x)}(u) = \left\{ 
\begin{matrix} 
0  & \text{if}~ x \in C_x,  u \in D(x)
 \\
\infty  & \text{otherwise}
\end{matrix} 
\right.$$
for each $(x,u) \in \mathbb{R}^n \times \mathbb{R}^m$.
Note that $\delta_S$ is an indicator function of the set $S$. 
Let $f:\mathbb{R}^n \times \mathbb{R}^m \to \bar{\mathbb{R}}$ be given by
$f(x,u):= \mathcal{L}(x,u) + \tilde{\delta}(x,u)$. 
Then, since $D$ is nonempty, $f$ is proper. The sublevel sets of $f$ for $\alpha \in \mathbb{R}$ are defined as
\begin{align} \label{eq:level-p1}
    \{(x &,u)  \in  \mathbb{R}^n \times \mathbb{R}^m \,:\, f(x,u)  \leq \alpha \} \\
& =\operatorname{Graph}(D) \cap  \{(x,u) \in \mathbb{R}^n \times \mathbb{R}^m \,:\, \mathcal{L}(x,u) \leq \alpha \}. \nonumber
\end{align}
\normalsize
Since $\mathcal{L}$ is lower semicontinuous, from \cite[Theorem $1.6$]{rockafellar2009variational} its sublevel sets are closed. Since $\operatorname{Graph}(D)$ is closed, then the sublevel sets of $f$ in \eqref{eq:level-p1} are closed. 
Thus, from \cite[Theorem $1.6$]{rockafellar2009variational} $f$ is lower semicontinuous.
Since $D$ has convex values and $\mathcal{L}$ is convex, then $f$ is convex. Furthermore, since $\mathcal{L}$ is strictly convex in $u$, $f$ is strictly convex in $u$.
If either $D$ has bounded values or the sublevel sets of $\mathcal{L}$ in $u$ are bounded, then the level sets of $f$ in $u$ are bounded. Namely, let $\mathcal{B} \subset C_x$ be a bounded set and let
$L=\{(x,u)\,:\, x \in \mathcal{B}, \, u \in D(x)\},$
for $\alpha \in \mathbb{R}$. The sublevel sets of $f$ in $u$ are defined by $ \{(x,u) \in \mathcal{B} \times \mathbb{R}^m \,:\, f(x,u)  \leq \alpha \} =
 L \cap  \{(x,u) \in \mathcal{B} \times \mathbb{R}^m \,:\, \mathcal{L}(x,u) \leq \alpha \}$.
Using 
\ifitsdraft
Corollary \ref{cor:num2} in the Appendix 
\else 
\cite[Corollary $3.27$]{rockafellar2009variational}
\fi
and the fact that $f$ is sublevel bounded in $u$, we conclude that $f$ is level coercive in $u$. Then, using \cite[Theorem $3.26$]{rockafellar2009variational}, we conclude that $f^\infty(0,u)>0$ (see the horizon function in \cite[Definition $3.17$]{rockafellar2009variational}) for all $u\not=0$. Then, 
\ifitsdraft
using Corollary \ref{cor:num1} in the Appendix, 
\else 
using \cite[Corollary $7.43$]{rockafellar2009variational},
\fi
 we conclude that $\kappa^*$ is single valued on $\dom \kappa^*$ and it is continuous on the interior of its domain.
 Since $f$ is proper, lower semicontinuous, and sublevel bounded in $u$, using \cite[Theorem $1.9$]{rockafellar2009variational} we conclude that $\dom \kappa^* = C_x$.
\end{proof}
\vspace{-5 mm}
\section{Guaranteeing Continuity of the Optimal Safe Control Law} \label{sec:cont_opt_safe}

Building from the results in Sections \ref{sec:suff_con} and \ref{sec:cont_optimal}, we formulate conditions for synthesizing an optimal, safe, and continuous control law. The constraint map $D$ is induced using CBF (safety constraint) as well as control and state constraints.
\vspace{-3 mm}
\subsection{Continuous, Safe, and Optimal Control Law }
As indicated in Theorem \ref{thm1}, to ensure safety, the selection of the control law in the outer neighborhood of the zero sublevel set of the CBF should be restricted appropriately. 
Let the set-valued map $D:C_x \rightrightarrows C_u$ be defined as  
\begin{align} 
D(x) := \left\{ 
\begin{matrix}
\bar{D}_{0}(x)  & \text{if}~ x \in S_1
 \\
 \mathcal{U}(x)  & \text{otherwise},
\end{matrix} 
\right.
\end{align}
where $\bar{D}_0 \subset \cal U$ is given in the proof of Theorem \ref{thm1} and indicates all safe feedback control laws on the set $S_1$, and $S_1:=\cl(U(K)\backslash K \cap C_x)$ indicates the corresponding outer neighborhood of the set $K$. Then, the set-valued map $D$ contains all feasible and safe feedback control laws. 
To select a continuous control law from $D$ by minimizing a cost function, based on Lemma  \ref{lm:opt_select}, $D$ should be continuous.

As explained in Remark \ref{rm:ncon_D}, the map $\bar{D}_{0}$ generally is not continuous. 
Here, we construct a continuous set-valued map $\widetilde{\cal U}$ from $D$.
To do this, we should design $\widetilde{\cal U}$ such that it is a subset of $\bar{D}_{0}$ when $ x \in S_1$ and also it blends smoothly with the set-valued map $\cal U$ on the boundary of $S$.

\begin{assumption} \label{ass:U2}
The set-valued map $\mathcal{U}$ is outer semicontinuous.
\end{assumption}

In the following result, we formulate conditions guaranteeing the selection of a continuous control law by minimizing a cost function and simultaneously ensuring safety.

\begin{theorem} \label{lm:safeOptCont}
Consider $\Sigma_{\mathcal{U}}$ in \eqref{eq:sys_u} such that Assumptions \ref{ass:C}-\ref{ass:U2} hold.
Let $B : \mathbb{R}^n \rightarrow \mathbb{R}$ be a locally Lipschitz CBF with respect to $(X_o, X_u) \subset \mathbb{R}^n \times \mathbb{R}^n$ defining the set $K$ in \eqref{eq:k}. Suppose
\begin{enumerate}
    \item  \label{th3:item1}The function $g:C_x \times \C_u \to \mathbb{R}$ defined in \eqref{eq:g_ll} is convex in $u$,
    \item \label{th3:item2}There exists a neighborhood of $\partial K$ such that the set-valued map $D_0:C_x \rightrightarrows C_u$  defined in \eqref{eq:D-zero} for $\gamma$ identically zero is nonempty.

\NoIndent{Let $S_1 =\cl(U( K)\backslash K \cap C_x )$ be a neighborhood of $\partial K$ such that $D_0$ is nonempty and \eqref{eq:cbf_llip} holds.  
Suppose $\widetilde{\mathcal{U}}:C_x\rightrightarrows C_u$ satisfies the following properties:}
\item \label{th3:item4}$\widetilde{\mathcal{U}}$ is continuous with nonempty, closed, and convex values,

\item \label{th3:item5}For each $x \in S_1$, $      \widetilde{\mathcal{U}}(x) \subset \bar{D}_0(x)$,

\item \label{th3:item3}For each $x \in C_x$, $\widetilde{\mathcal{U}}(x) \subset \mathcal{U}(x)$.
\end{enumerate}
Let $\mathcal{L}:\mathbb{R}^n \times \mathbb{R}^m \to \mathbb{R}$ be proper, lower semicontinuous, convex in both arguments, and strictly convex in its second argument.
Let $\kappa^*:C_x \to C_u$ be defined by
\begin{equation} \label{eq:opt}
\begin{aligned}
  \kappa^*(x)=\argmin_{} \quad &  \mathcal{L}(x,u)\\
\textrm{s.t.} \quad  u \in \widetilde{\mathcal{U}}(x). \\
\end{aligned}
\end{equation}
If one of the following conditions holds
\begin{enumerate}
    \item $\mathcal{L}$ is sublevel bounded in $u$;
    \item The set-valued map $\widetilde{\mathcal{U}}$ has bounded values,
\end{enumerate}
then $\kappa^*$ is continuous on $\operatorname{int}(C_x)$, and the resulting differential inclusion $\Sigma_{\kappa^*}$ is safe with respect to $(X_o, X_u)$, namely, $\kappa^*$ is the optimal safe control law.
\end{theorem}
\begin{proof}
 Lemma \ref{lm:opt_select} implies that $\kappa^*$ is single valued on $C_x$ and continuous on $\operatorname{int}(C_x).$
Since $\kappa^*$ is a continuous selection from a subset of $\bar{D}_0$ on $U(\partial K)\backslash K$, then $\Sigma_{\kappa^*}$ satisfies
\ifitsdraft
\begin{align*}
    \langle \zeta, \eta\rangle \leq 0 \quad \forall x \in \left(U(\partial K)\backslash K \right)\cap C_x,& \,\zeta \in \partial_C B(x),\\
    &\,\forall \eta \in F(x,\kappa^*(x)).
\end{align*}
\else
$ \langle \zeta, \eta\rangle \leq 0$ for each $x \in \left(U(\partial K)\backslash K \right)\cap C_x, \,\zeta \in \partial_C B(x)$, and for each $ \eta \in F(x,\kappa^*(x)).$
\fi
Therefore, \cite[Theorem $4$]{maghenem2021sufficient} implies that $K$ is forward pre-invariant for $\Sigma_{\kappa^*}$. Then, 
$\Sigma_{\kappa^*}$ is safe w.r.t. $(X_o, X_u)$.
\end{proof}
 \vspace{-5 mm}
\subsection{A Sample Construction of \texorpdfstring{$\widetilde{\mathcal{U}}$}{TEXT}}

Because $D$ consists of two set-valued maps $\bar{D}_0$ and $\mathcal{U}$, to construct the continuous set-valued map $\widetilde{\mathcal{U}}$ from $D$, we need to apply two types of blending. 
In general, $\bar{D}_0$  is not continuous, as a result of discontinuity in $g$ with respect to $x$. In obstacle avoidance, for example,  when the barrier function is defined as the minimum or maximum of some hyperplanes, there are discontinuities in $g$ as different constraints are active in the different regions around the obstacle. 
First, we should find some continuous map $D_0^s:S_1 \to C_u$ such that $D_0^s(x) \subset \bar{D}_0(x)$ for each $x \in S$.
Second,  we must blend continuously two continuous maps, $D_0^s$ and $\cal U$.

In Example \ref{ex:num1}, we present an approach to smoothen $\bar{D}_0$ when the unsafe set is defined as a system of linear inequalities, and the barrier function is defined as the minimum of the hyperplanes corresponding to the unsafe set.

Here, we present an approach for continuously blending two continuous maps using the Minkowski sum. Let $S \subset \reals^n$ and $\epsilon: S \to \mathbb{R}_{>0}$ be continuous. The $\epsilon$-neighborhood of $S$, $U_{\epsilon} (S)$, is defined as
\ifitsdraft
\begin{align}
    U_{\epsilon} ( S):= \bigcup_{x \in S} (x + \epsilon(x) \mathbb{B}).
\end{align}
\else
$ U_{\epsilon} ( S):= \bigcup_{x \in S} (x + \epsilon(x) \mathbb{B})$.
\fi
\ifitsdraft
Given $S_1, S_2 \subset \reals^n$ and $\alpha_1,\alpha_2 \in \reals$, based on Minkowski sum of sets, define $\alpha_1 S_1 + \alpha_2 S_2:= \{\alpha_1 s_1+\alpha_2 s_2\,:\, s_1 \in S_1, \, s_2 \in S_2\}$.
\else 
\fi
\begin{lemma} \label{lm:lm5}
Consider a closed set $K \subset \reals^n$, and set-valued maps $F_1:K \rightrightarrows \mathbb{R}^m$ and $F_2:U_{\epsilon} (\partial K) \rightrightarrows \mathbb{R}^m$ that are continuous with closed and convex values, where  $\epsilon:\partial K \to \reals_{>0}$ is a continuous function such that $\epsilon(x)> \epsilon_1$ for each $x \in \partial K$, where $\epsilon_1 \in \mathbb{R}_{>0}$. 
Then, the set-valued map $G:K \cup U_{\epsilon} (\partial K) \rightrightarrows \reals^m$ defined as 
\begin{align} \label{eq:GG}
    G(x):=\left\{ 
\begin{matrix} 
F_1(x)  & \text{if}~ x \in K \backslash U_{\epsilon_1} (\partial K)
 \\
 F_3(x) & \text{if}~ x \in K \cap U_{\epsilon_1} (\partial K)
 \\
 F_2(x)  & \text{if}~ x \in U_{\epsilon} (\partial K) \backslash K
\end{matrix} 
\right.
\end{align}
where $F_3(x):=   \frac{d(x,\partial K)}{\epsilon_1} F_1(x) + (1-\frac{d(x,\partial K)}{\epsilon_1}) F_2(x)$, is lower and outer semicontinuous with closed and convex values.
\end{lemma}
\begin{proof}
 Since $\partial K$ is a closed set, then the distance function $x \to d(x,\partial K)$ is continuous \cite{rockafellar2009variational}. Therefore, $F_3$ is continuous. Since $F_1$ and $F_2$ have closed and convex values, using the Minkowski sum of two closed and convex sets is closed and convex \cite{schneider2014convex}, we conclude that $F_3$ has closed and convex values. For each $x \in \partial K$, $F_3(x) = F_2(x)$ and for each $x \in \partial (K\backslash U_{\epsilon_1}(\partial K))$, $F_3(x) = F_1(x)$, therefore, $G$ is continuous with closed and convex values.
\end{proof}
\blue{\begin{example} \label{ex:num1}
Consider the system
\ifitsdraft
\begin{equation*}
   \Sigma_{\mathcal{U}}: \quad \dot{x} =  Ax+Bu \quad   x \in \mathbb{R}^2.
\end{equation*}
\else
$\Sigma_{\mathcal{U}}: \dot{x} = Ax+Bu$ for $x \in \mathbb{R}^2$.
\fi
Suppose we have a rectangular shape obstacle.
Let $p_0$ be the central point of the obstacle and, for each $i\in \{1,\cdots,4\}$, $q_i$  denotes the middle point of each edge of the obstacle. The unsafe set is defined as the intersection of the halfspaces 
$(p_0 - q_i)^{\top}(x-q_i) > 0$ for $i=1,\cdots,4$. We have
\ifitsdraft
\begin{align}
    X_u &=\{x \in \mathbb{R}^2 \,:\, A_u x  > b_u \}
\end{align}
where $A_u=\begin{pmatrix} (p_0 - q_1)^{\top} \\(p_0 - q_2)^{\top}\\(p_0 - q_3)^{\top}\\(p_0 - q_4)^{\top}
\end{pmatrix}$, and $b_u=\begin{pmatrix} (p_0 - q_1)^{\top} q_1\\(p_0 - q_2)^{\top} q_2\\(p_0 - q_3)^{\top} q_3\\(p_0 - q_4)^{\top} q_4
\end{pmatrix}$.
\else
$X_u =\{x \in \mathbb{R}^2 \,:\, A_u x  > b_u \}$,
where $A_u=\big(p_0 - q_1,p_0 - q_2,p_0 - q_3,p_0 - q_4 \big)^{\top}$, and $b_u=\big((p_0 - q_1)^{\top}q_1,(p_0 - q_2)^{\top}q_2,(p_0 - q_3)^{\top}q_3,(p_0 - q_4)^{\top}q_4 \big)^{\top}$. 
\fi
Let $\mathcal{U}:\mathbb{R}^2 \rightrightarrows \mathbb{R}^2$ be defined as $\mathcal{U}(x) =  [-5,5]^2 $.

We set the initial set, $X_o$, to be the complement of the unsafe set with an extra distance as $X_o =  \mathbb{R}^2 \backslash \{x \in \mathbb{R}^2 \,:\, A_u x  > b_u -d \mathds{1}\}$,
where $\mathds{1}$ denotes a vector of ones, and $d$ denotes the extra distance.
\ifitsdraft
Next, we define $K=X_0$ and CBF as 
\begin{equation*}
    B(x):= \min_{i \in \{1,\cdots,4\}} B_i(x) \quad \forall x \in \mathbb{R}^2,
\end{equation*}
where the $B_i$ is given by
\begin{equation*}
    B_i(x):= (p_0 - q_i)^{\top}(x-q_i)+d \quad i \in \{1,\cdots,4\}.
\end{equation*}
The gradient of the $B_i$ is
\begin{equation} \label{eq:gb}
    \nabla B_i(x)=  (p_0 - q_i) \quad i \in \{1,\cdots,4\}.
\end{equation}
Based on \cite[Proposition $2.3.12$]{clarke1990optimization}, the Clarke generalized gradient of $B$ is given by
\begin{equation*}
    \partial_C B(x)=  \operatorname {co} (\{\nabla B_i(x)\,:\, i \in \mathcal{I}(x)\})
\end{equation*}
where $\mathcal{I}:\mathbb{R}^2 \rightrightarrows \{1,\cdots,4\}$ indicates the active $B_i$'s 
$$\mathcal{I}(x) :=\{i \,:\, B(x)=B_i(x), i \in \{1,\cdots,4\}\}.$$
Therefore, $\partial_C B(x)$ is given by
\begin{align*}
    \partial_C B(x) &= \{\sum_{i \in \mathcal{I}(x)} \theta_i \nabla B_i(x)\,:\,  \theta_i \geq 0,\,\sum_{i \in \mathcal{I}(x)} \theta_i = 1\}
\end{align*}
\else
Next, we define $K=X_0$ and the CBF as $B(x):= \min_{i \in \{1,\cdots,4\}} B_i(x)$ for each $x \in \mathbb{R}^2$, where the $B_i$ is given by
$B_i(x):= (p_0 - q_i)^{\top}(x-q_i)+d$ for $i \in \{1,\cdots,4\}$.
Based on \cite[Proposition $2.3.12$]{clarke1990optimization}, the Clarke generalized gradient of $B$ is given by
$\partial_C B(x) = \{\sum_{i \in \mathcal{I}(x)} \theta_i \nabla B_i(x)\,:\,  \theta_i \geq 0,\,\sum_{i \in \mathcal{I}(x)} \theta_i = 1\}$ where 
 $\mathcal{I}:\mathbb{R}^2 \rightrightarrows \{1,\cdots,4\}$ indicates the active $B_i$s and is defined as 
$\mathcal{I}(x) :=\{i \,:\, B(x)=B_i(x), i \in \{1,\cdots,4\}\}$.
\fi
The function $g$ from Theorem \ref{lm:safeOptCont} is
\ifitsdraft
\begin{equation}\label{eq:g_opt}
\begin{aligned} 
   & g(x,u)=\sup_{ \zeta \in \partial_C B(x)} \langle \zeta, Ax+Bu\rangle
     \\
    &=\sup_{} \quad   \langle\sum_{i \in \mathcal{I}(x)} \theta_i (p_0 - q_i),Ax+Bu\rangle   \\
& \quad \qquad\textrm{s.t.} \quad  \theta_i \geq 0 \quad i \in \mathcal{I}(x),
   \,\,\sum_{i \in \mathcal{I}(x)} \theta_i = 1
\end{aligned}
\end{equation}
The optimization \eqref{eq:g_opt} is over $\theta$ and of the form of a linear function over the probability simplex. 
\else
$g(x,u)=
    \sup_{} \{   \langle \sum_{i \in \mathcal{I}(x)} \theta_i (p_0 - q_i),Ax+Bu \rangle  \,:\,
   \theta_i \geq 0, \,\, i \in \mathcal{I}(x),\,\,
   \sum_{i \in \mathcal{I}(x)} \theta_i = 1 \}.$
The above optimization is over $\theta$ and of the form of a linear function over the probability simplex.
\fi
\ifitsdraft
Then, the dual problem for \eqref{eq:g_opt} is 
\begin{equation}
\begin{aligned}
   & \min_{} \quad   \nu  \\
& \quad  \textrm{s.t.} \quad  \nu \geq \langle  p_0-q_i,Ax+Bu\rangle \quad \forall i \in \mathcal{I}(x).
\end{aligned}
\end{equation}
Therefore $\nu = \max_{i \in \mathcal{I}(x)} \langle  p_0-q_i,u\rangle$ and the optimal solution to the primal problem is $\theta_{i^*} =1$, where $i^*$ is the index of the maximum $\langle  p_0-q_i,Ax+Bu\rangle$ for $i=\{1,\cdots,4\}$.
Therefore,
\begin{align*}
    g(x,u) &= \max \{\langle \nabla B_i(x), Ax+Bu \rangle, i \in \mathcal{I}(x)\}.
\end{align*}
\else
Therefore, $g(x,u) = \max \{\langle \nabla B_i(x), Ax+u \rangle, i \in \mathcal{I}(x)\}$.
\fi
Since $g$ is the pointwise maximum of affine functions in $u$, it is convex in $u$.
Then, $D_0$ in Theorem \ref{lm:safeOptCont}  is given by
\ifitsdraft
\begin{align*}
    D_0(x) &= \{u \in \mathcal{U}(x) \,:\, g(x,u)<0\} \\
    &=\{u \in [-5,5]^2 \,:\, \langle p_0-q_i, Ax+Bu \rangle<0,\, \forall i \in \mathcal{I}(x)\}.
\end{align*}
We observe that since, based on the position of system w.r.t. the obstacle, one or two of the constraints are active, the map $D_0$ is nonempty.
\else
$D_0(x) =\{u \in [-5,5]^2 \,:\, \langle p_0-q_i, Ax+Bu \rangle<0,\, \forall i \in \mathcal{I}(x)\}$. 
We observe that the map $D_0$ is nonempty.
\fi
First, we want to smoothen the discontinuities of  $\bar{D}_0$, which is induced by changing the active constraints.  
We define $\mathcal{I}_{\alpha}$ for some $\alpha >0$ as 
$\mathcal{I}_{\alpha}(x) :=\{i \,:\, |B(x)-B_i(x)| \leq \alpha, i \in \{1,\cdots,4\}\}$.
 When $\mathcal{I}(x)=\mathcal{I}_{\alpha}(x)$, the conditions to satisfy are
 $ \nabla B_i(x)^{\top} (Ax+Bu) \leq 0$, for each $i \in \mathcal{I}(x)$. Let $\phi(x):= \frac{\pi |B(x)-B_j(x)|}{2 \alpha}$ for $j \in \mathcal{I}_\alpha(x) \backslash \mathcal{I}(x)$. Then, 
\ifitsdraft
we have 
\begin{align*}
    \big( \sin(\phi)\nabla B_i(x)  + \cos(\phi)  \nabla & B_j(x)  \big)^{\top}(Ax+Bu) \leq 0, 
    \\
    &\quad \forall i \in \mathcal{I}(x), j \in \mathcal{I}_\alpha(x) \backslash \mathcal{I}(x)
\end{align*}
\else
for each $i \in \mathcal{I}(x), j \in \mathcal{I}_\alpha(x) \backslash \mathcal{I}(x)$ we have $\big( \sin(\phi(x))\nabla B_i(x)  + \cos(\phi(x))  \nabla  B_j(x)  \big)^{\top}(Ax+Bu) \leq 0$.
\fi
At $x$'s such that $\mathcal{I}$ is not equal to $\mathcal{I}_\alpha$, means that we are near to changing the active constraints. Using $\sin$ and $\cos$ functions, we smoothen these transitions.  
We define $\widetilde{\mathcal{U}}(x) = 
    \{u\in [-5,5]^2 \,:\, H(x)^{\top} u \leq h(x)  \},$
where $H$ is defined based on the conditions we explained above as
\ifitsdraft
\begin{equation}
    H(x) = \begin{pmatrix}   
    (p_0-q_i)^{\top}B \\
    (\sin(\phi)(p_0-q_i)+\cos(\phi)(p_0-q_j))^{\top}B
    \end{pmatrix}
\end{equation}
where $\forall i \in \mathcal{I}(x), j \in \mathcal{I}_\alpha(x)\backslash \mathcal{I}(x)$.
The $h$ is defined as
\begin{equation}
    h(x) = \begin{pmatrix}   
    -(p_0-q_i)^{\top}Ax  \\
    -(\sin(\phi)(p_0-q_i)+\cos(\phi)(p_0-q_j))^{\top}Ax 
    \end{pmatrix}
\end{equation},
\else
$H(x) = \begin{pmatrix}   
    (p_0-q_i)^{\top}B \\
    (\sin(\phi)(p_0-q_i)+\cos(\phi)(p_0-q_j))^{\top}B
    \end{pmatrix}$
for each $ i \in \mathcal{I}(x)$, and $ j \in \mathcal{I}_\alpha(x)\backslash \mathcal{I}(x)$.
The $h$ is defined as
$h(x) = \begin{pmatrix}   
    -(p_0-q_i)^{\top}Ax  \\
    -(\sin(\phi)(p_0-q_i)+\cos(\phi)(p_0-q_j))^{\top}Ax 
    \end{pmatrix}$,
\fi
and we add $\max\{0,-M B(x)\}$ to every row of $h$ for some sufficiently large positive number $M$ to make the constraint $H(x)^{\top} u \leq h(x)$ ineffective when we are in the safe region. We used this method instead of the one in Lemma \ref{lm:lm5} since that can result in a set-valued map with nonconvex values.
 Considering that all the constraints defining $\widetilde{\mathcal{U}}$ are affine with respect to $u$ at each $x$, and the constraints with strict inequality are nonempty, \cite[Example $5.10$]{rockafellar2009variational} implies that  $\widetilde{\mathcal{U}}$ is continuous.
Furthermore, $\widetilde{\mathcal{U}}$ has nonempty and compact values. 
Theorem  \ref{lm:safeOptCont} implies that if the cost function $\mathcal{L}$ is lower semicontinuous, convex and strictly convex in $u$, then the optimal control law $\kappa^{*}$ is continuous and safe. 
Figure \ref{fig:traj2} shows trajectories using 
$\mathcal{L}(x,u) = \frac{1}{2}u^{\top} u$ and control Lyapunov function $V(x):=\frac{1}{2}x^2$ for $A = \begin{pmatrix}0&1\\-1&-1\end{pmatrix}$, $B = I$,
$\alpha = 0.01$, and $M=100$.
\ifitsdraft
\begin{figure}[t]
\centering
\includegraphics[width=3.5in, height=2.5in, trim = {0 1.6cm 0 1.5cm}]{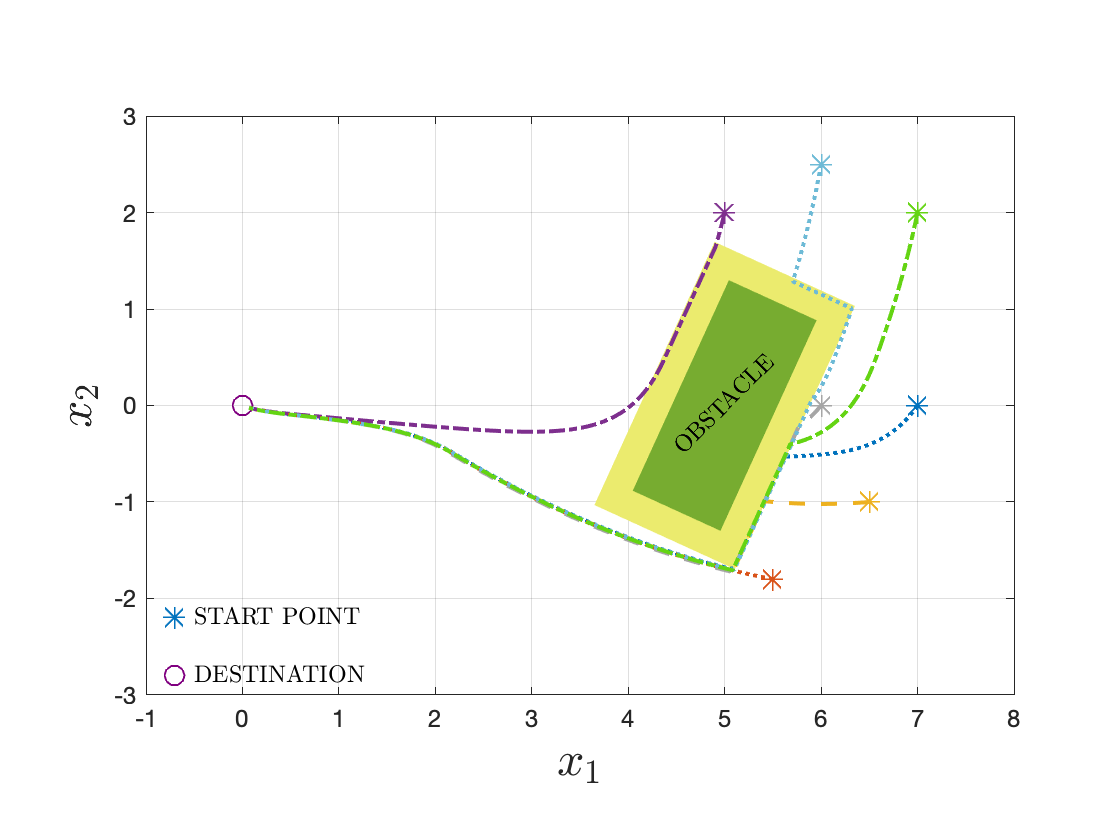}
\caption{Trajectories for Example \ref{ex:num1}. The obstacle is the green rectangle, and the initial set is the area outside of the yellow and green regions.}
\label{fig:traj2}
\end{figure} 
\else
\begin{figure}[t]
\centering
\includegraphics[width=3.0in, height=2.0in, trim = {0 1.6cm 0 1.5cm}]{lin_traj.png}
\caption{The obstacle is the green rectangle, and the initial set is the area outside of the yellow and green regions.}
\label{fig:traj2}
\end{figure} 
\fi
\end{example}}
\vspace{-6 mm}
\section{Conclusion} \label{sec:conlusion}

This paper studies controlled safety of constrained differential inclusions using nonsmooth CBFs. We develop sufficient conditions to select a continuous control law using CBFs. 
Furthermore, we study conditions to find optimal safe control laws while minimizing the cost function. We illustrate the results in an obstacle avoidance example.
\ifitsdraft
We extend the results for hybrid control systems and consider conditions for robust controlled safety for future work. 
\else
\fi

\ifitsdraft
\appendix
\begin{lemma} \label{lm:michael_extend}
(Selection Theorem \cite{michael1956continuous}) 
Consider a closed set $K \subset \mathbb{R}^n$ and a set-valued map $\Phi : K \rightrightarrows \mathbb{R}^m$ that is lower semicontinuous with $\Phi(x)$ nonempty, closed, and convex for all $x \in K$. Then, for each closed set $A \subset K$ and for each $\phi : A \rightarrow \mathbb{R}^m$ that is a selection of $\Phi$ on $A$; namely, $\phi$ is continuous and 
$ \phi(x) \in \Phi(x)$ for every $ x \in A$,
there exists $\phi_e : K \rightarrow \mathbb{R}^m$ extending 
$\phi$ to $K$ that is a selection of $\Phi$ on $K$; namely, 
$\phi_e$ is continuous,  
$ \phi_e(x) = \phi(x)$ for every $ x \in A, $
and
$ \phi_e(x) \in \Phi(x)$ for every $ x \in K. $ 
\end{lemma}

\begin{theorem} \label{th:att}
(Attainment of a minimum \cite[Theorem $1.9$]{rockafellar2009variational})
Suppose $f:\mathbb{R}^n \to \mathbb{R}$ is lower semicontinuous,
sublevel bounded and proper. Then the value $\inf f$ is finite and the
set $\argmin f$ is nonempty and compact.
\end{theorem}
\begin{lemma} \label{lm:appOSC_usc}
\cite[Lemma $5.15$]{goebel2012hybrid} Let $F:\mathbb{R}^n \rightrightarrows \mathbb{R}^m$ be a upper semicontinuous set-valued mapping. Consider $x \in \mathbb{R}^n$ such that $F(x)$ is closed. Then, $F$ is outer semicontinuous at $x$. If $F$ is locally bounded as $x$, then the reverse implication is true.
\end{lemma}

\begin{lemma} \label{lm:mul_sing}
\cite[Corollary $5.20$]{rockafellar2009variational} For any single-valued mapping $F:\mathbb{R}^n \to \mathbb{R}^m$, viewed as a special case of a set-valued mapping, the following properties are equivalent:
\begin{itemize}
    \item $F$ is continuous at $x$;
    \item $F$ is outer semicontinuous and locally bounded at $x$;
    \item $F$ is inner semicontinuous at $x$.
\end{itemize}
\end{lemma}

\begin{corollary} \label{cor:num1}
\cite[Corollary $7.43$]{rockafellar2009variational}
Let $f:\mathbb{R}^n \times \mathbb{R}^n \to \mathbb{R}\cup \{\pm \infty\}$ be proper, lower semicontinuous, convex, and such that $f^\infty(0,u)>0$ for all $u \not= 0$. Suppose 
$\kappa(x)=\argmin_u f(x,u).$
If $f(x,u)$ is strictly convex in $u$, then $\kappa$ is single valued on $\dom \kappa$ and continuous on $\operatorname{int}(\dom \kappa).$
\end{corollary}

\begin{corollary} \label{cor:num2}
\cite[Corollary $3.27$]{rockafellar2009variational}
For any proper, lower semicontinuous, and convex function $f$ on $\mathbb{R}^n$, level coercivity is equivalent level boundedness.
\end{corollary}
\else
\fi

\bibliography{biblio}
\bibliographystyle{ieeetr}     
\end{document}